\documentclass[11pt, a4paper, reqno]{amsart}
\usepackage[colorlinks,bookmarksopen=false]{hyperref}
\usepackage[english]{babel}
\usepackage{amsmath}
\usepackage{amssymb}
\usepackage{amsthm}
\usepackage{latexsym}
\newtheorem{prop}{Proposition}
\newtheorem{thm}[prop]{Theorem}
\newtheorem{cor}[prop]{Corollary}

\theoremstyle{definition}
\newtheorem*{defn*}{Definition}

\newtheorem*{remark*}{Remark}
\newtheorem*{remarks*}{Remarks}
\newtheorem{remark}[prop]{Remark}
\newtheorem*{scholium*}{Scholium}
\numberwithin{equation}{section}

\newcommand{\RR}{\mathbf{R}}
\newcommand{\ZZ}{\mathbf{Z}}

\begin{document}
\title[Strong law of large numbers with concave moments]{Strong law of large numbers\\ with concave moments}
\author{Anders Karlsson}
\address{Royal Institute of Technology, Sweden / EPFL and Geneva, Switzerland}
\email{akarl@math.kth.se / nicolas.monod@epfl.ch}
\author{Nicolas Monod}
%
\begin{abstract}
It is observed that a wellnigh trivial application of the ergodic theorem from~\cite{Karlsson-Ledrappier06} yields a strong LLN for arbitrary concave moments.
\end{abstract}
\maketitle
\let\languagename\relax  
\thispagestyle{plain}
\pagestyle{plain}
\vskip-7mm
\noindent
\frame{\parbox[t]{\textwidth}{
\bfseries
Not for publication: we found that Aaronson--Weiss essentially proved Theorem~\ref{thm:main}, see J.~Aaronson, \emph{An introduction to infinite ergodic theory} (AMS Math. Surv. Mon.~50, 1997), pages~65--66.
}}
\section{Introduction}

Let $\Omega$ be a standard probability space and $L:\Omega \to\Omega$ and ergodic measure-preserving transformation. Let $f:\Omega\to\RR$ be any measurable map and consider the \emph{Birkhoff sums} $S_n=\sum_{k=0}^{n-1} f\circ L^k$. Recall that this ergodic (or ``stationary'') setting includes the case of the sums $\sum_{k=0}^{n-1} X_k$ of a family $\{X_k\}_k$ of i.i.d.\ random variables.

\begin{thm}\label{thm:main}
Let $D:\RR_+\to\RR_+$ be any concave function with $D'(\infty)=0$.\\ If $D(|f|)$ is integrable, then
$$~\kern20mm\lim_{n\to\infty} \frac{1}{n} D(|S_n|)\ =\  0 \kern7mm \text{a.s.}$$
\end{thm}

\begin{remark}\label{rem:main}
(i)~For the notation $D'(\infty)$, recall that the derivative $D'$ exists except possibly on a countable set and is non-increasing. Moreover, $\lim\limits_{t\to\infty}D(t)/t= D'(\infty)$.\\
(ii)~The condition $D'(\infty)=0$ is not a restriction, since otherwise straightforward estimates reduce the question to Birkhoff's theorem and $D(|S_n|)/n$ tends to $D'(\infty)|\int\! f|$.
\end{remark}

Playing with the choice of the arbitrary concave function $D$, one gets examples old and new. For instance, $D(t)=t^p$ yields a Marcinkiewicz--Zygmund theorem:

\begin{cor}
Let $0<p<1$. If $f\in L^p$, then $\displaystyle\lim_{n\to\infty} \frac{1}{n^{1/p}} S_n=0$.\hfill\qedsymbol
\end{cor}

\begin{remark}
Marcinkiewicz--Zygmund~\cite[\S6]{Marcinkiewicz-Zygmund37} work in the i.i.d.\ case but with $0<p<2$; however, when $p>1$, no such statement can hold in the ergodic generality even under the strongest assumptions, see Proposition~2.2 in~\cite{Baxter-Jones-Lin-Olsen}. The independence condition was removed by S.~Sawyer~\cite[p.~165]{Sawyer66}. The beautiful geometric proof by Ledrappier--Lim~\cite{Ledrappier-Lim} for $p=1/2$ has inspired the present note.
\end{remark}

Recall that $f$ is \emph{log-integrable} if $\log^+|f|\in L^1$, where $\log^+=\max(\log, 0)$. The choice $D(t)=\log(1+t)$ yields:

\begin{cor}\label{cor:log_p}
If $f$ is log-integrable, then $\displaystyle\lim_{n\to\infty}  |S_n|^{1/n}=1$.\hfill\qedsymbol
\end{cor}

Observe that for functions $\RR_+\to\RR_+$, concavity is preserved under composition. Combining this operation and shift of variables, one has further examples such as:

\begin{cor}
Let $p>0$. If $\log^+|f|\in L^p$, then $\displaystyle\lim_{n\to\infty}  |S_n|^{n^{-1/p}}=1$.
\end{cor}

\section{Proofs}
\emph{Theorem~\ref{thm:main}}. Notice that $D$ is non-decreasing and subadditive. We can assume $D(0)=0$ upon adding a constant. We can assume that $D$ tends to infinity since otherwise it is bounded. Thus, $D(|x-y|)$ defines a proper invariant metric on the group $\RR$ and $S_n$ is the random walk associated to $(\Omega, L)$. Recall that a \emph{horofunction} $h$ (normalized by $h(0)=0$) is any limit point in the topology of compact convergence of the family $D(|t-x|) - D(|t|)$ of functions of $x$ indexed by $t$. According to~\cite{Karlsson-Ledrappier06}, there are horofunctions $h_\omega$ on $\RR$ such that
$$\lim_{n\to\infty} \frac{1}{n} D(S_n(\omega))\ =\  - \lim_{n\to\infty} \frac{1}{n} h_\omega(S_n(\omega))\kern7mm \mathrm{a.s.}(\omega)$$
However, $D'(\infty)=0$ implies that $h=0$ is the only horofunction since
$$\lim_{t\to\pm\infty} \big(D(|t-x|) - D(|t|)\big)\ = 0\kern7mm \forall\,x. \eqno{\text\qedsymbol}$$

\smallskip

\emph{Corollary~\ref{cor:log_p}}. For any $p>0$ there is $x_0$ such that $D(t)=(\log(x+x_0))^p$ is concave on $\RR_+$. Now $D(|f|)$ is integrable and the statement follows.\hfill\qedsymbol

\section{Comments and references}
(i)~We used only a very special case of the LLN of~\cite{Karlsson-Ledrappier06}, which applies to group-valued random variables. Theorem~\ref{thm:main} holds indeed also in that setting with identical proof, but can immediately be reduced to the real-valued case.

\smallskip

(ii)~The point of the present note is that the LLN of~\cite{Karlsson-Ledrappier06} brings new insights even when the group is the range $\RR$ of classical random variables, since we can endow it with various invariant metrics. There is indeed a wealth of such metrics; recall that even $\ZZ$ admits an invariant metric whose completion is Urysohn's universal polish space~\cite{Urysohn27} (by Theorem~4 in~\cite{Cameron-Vershik}). Wild \emph{proper} metrics can be constructed by means of weighted infinite generating sets.

\smallskip

(iii)~One can relax the concavity assumption is various ways. For instance, keeping $D(t)/t\to 0$, it suffices to assume that $D$ is quasi-concave in the sense that Jensen's inequalities hold up to a multiplicative constant. Indeed, this implies that $D$ can be constrained within two proportional concave functions~\cite[Theorem~1]{Mulholland32}.

\smallskip

(iv)~V.~Petrov~\cite{Petrov96} discusses laws of the form $S_n/a_n\to 0$ in the i.i.d.\ case.


{\small
\begin{thebibliography}{1}

\bibitem{Baxter-Jones-Lin-Olsen}
J.~Baxter, R.~Jones, M.~Lin, and J.~Olsen.
\newblock S{LLN} for weighted independent identically distributed random
  variables.
\newblock {\em J. Theoret. Probab.}, 17(1):165--181, 2004.

\bibitem{Cameron-Vershik}
P.~J. Cameron and A.~M. Vershik.
\newblock Some isometry groups of the {U}rysohn space.
\newblock {\em Ann. Pure Appl. Logic}, 143(1-3):70--78, 2006.

\bibitem{Karlsson-Ledrappier06}
A.~Karlsson and F.~Ledrappier.
\newblock On laws of large numbers for random walks.
\newblock {\em Ann. Probab.}, 34(5):1693--1706, 2006.

\bibitem{Ledrappier-Lim}
F.~Ledrappier and S.~Lim.
\newblock A proof of a {$L^{1/2}$} ergodic theorem.
\newblock preprint.

\bibitem{Marcinkiewicz-Zygmund37}
J.~Marcinkiewicz and A.~Zygmund.
\newblock {Sur les fonctions ind\'ependantes}.
\newblock {\em Fundam. Math.}, 29:60--90, 1937.

\bibitem{Mulholland32}
H.~P. Mulholland.
\newblock {The generalization of certain inequality theorems involving powers.}
\newblock {\em Proc. Lond. Math. Soc., II. Ser.}, 33:481--516, 1932.

\bibitem{Petrov96}
V.~V. Petrov.
\newblock On the strong law of large numbers.
\newblock {\em Statist. Probab. Lett.}, 26(4):377--380, 1996.

\bibitem{Sawyer66}
S.~A. Sawyer.
\newblock Maximal inequalities of weak type.
\newblock {\em Ann. of Math. (2)}, 84:157--174, 1966.

\bibitem{Urysohn27}
P.~Urysohn.
\newblock {Sur un espace m\'etrique universel (r\'edig\'e par P.~Alexandroff)}.
\newblock {\em Bulletin Sc. Math. (2)}, 51:43--64 \& 74--90, 1927.

\end{thebibliography}
}
\end{document}